\documentclass[11pt]{amsart}
\usepackage{a4}

\title{Some properties of the pseudo-Smarandache function}

\author{Richard Pinch}
\address{2 Eldon Road, Cheltenham, Glos GL52 6TU, U.K.}
\email{rgep@chalcedon.demon.co.uk}
\date{2 April 2005}

\newcommand{\e}{\mbox{e}}
\renewcommand{\d}{{\rm d}}

\def\abs|#1|{{\left\vert{#1}\right\vert}}
\def\paren(#1){{\left({#1}\right)}}
\def\floor[#1]{{\left\lfloor{#1}\right\rfloor}}
\def\ceiling[#1]{{\left\lceil{#1}\right\rceil}}
\def\congruent{\equiv}

\def\divides{\mid}
\def\Larrow{\hbox to 3cm{\leftarrowfill}}
\def\Rarrow{\hbox to 3cm{\rightarrowfill}}

\newtheorem{theorem}{Theorem}
\newtheorem{lemma}{Lemma}

\begin{document}

\thispagestyle{empty}
\begin{abstract}
Charles Ashbacher \cite{Ash:pluckings} has posed a number of questions relating to
the pseudo-Smarandache function $Z(n)$.  In this note we show that the ratio of
consecutive values $Z(n+1)/Z(n)$ and $Z(n-1)/Z(n)$ are unbounded; that $Z(2n)/Z(n)$
is unbounded; that $n/Z(n)$ takes every integer value infinitely often; and that the
series $\sum_n 1/Z(n)^\alpha$ is convergent for any $\alpha > 1$.
\end{abstract}

\maketitle

\section{Introduction}

We define the $m$-th triangular number $T(m) = \frac{m(m+1)}{2}$.  Kashihara
\cite{Kash:smarandache} has defined the {\em pseudo-Smarandache} function $Z(n)$ by
$$
Z(n) = \min \{ m : n \divides T(m) \} .
$$

Charles Ashbacher \cite{Ash:pluckings} has posed a number of questions relating to
the pseudo-Smarandache function $Z(n)$.  In this note we show that the ratio of
consecutive values $Z(n)/Z(n-1)$ and $Z(n)/Z(n+1)$ are unbounded; that $Z(2n)/Z(n)$
is unbounded; and that $n/Z(n)$ takes every integer value infinitely often.
He notes that the series $\sum_n 1/Z(n)^\alpha$ is divergent for $\alpha=1$ and 
asks whether it is convergent for $\alpha=2$.  He further suggests that
the least value of $\alpha$ for which the series converges ``may never be known'' .
We resolve this problem by showing that the series converges for all $\alpha>1$.

\section{Some properties of the pseudo-Smarandache function}

We record some elementary properties of the function $Z$.

\begin{lemma}

\begin{enumerate}
\item If $n \ge T(m)$ then $Z(n) \ge m$.  $Z(T(m)) = m$.
\item For all $n$ we have $\sqrt{n} < Z(n)$.
\item $Z(n) \le 2n-1$, and if $n$ is odd then $Z(n) \le n-1$.
\item If $p$ is an odd prime dividing $n$ then $Z(n) \ge p-1$.
\item $Z(2^k) = 2^{k+1}-1$.
\item If $p$ is an odd prime then $Z(p^k) = p^k-1$ and
$Z(2p^k) = p^k-1$ or $p^k$ according as $p^k \congruent 1$ or $3 \bmod 4$.
\end{enumerate}
\end{lemma}

We shall make use of Dirichlet's Theorem on primes in arithmetic progression in the
following form.

\begin{lemma}
Let $a, b$ be coprime integers.  Then the arithmetic progression $a+bt$ is prime
for infinitely many values of $t$.
\end{lemma}

\section{Successive values of the pseudo-Smarandache function}

Using properties (3) and (5), Ashbacher observed that $\abs|Z(2^k)-Z(2^k-1)| > 2^k$
and so the difference between the conecutive values of $Z$ is unbounded.  He asks
about the ratio of consecutive values.

\begin{theorem}
For any given $L > 0$ there are infinitely many values of $n$ such that $Z(n+1)/Z(n) > L$,
and there are infinitely many values of $n$ such that $Z(n-1)/Z(n) > L$.
\end{theorem}
\begin{proof}
Choose $k \congruent 3 \bmod 4$, so that $T(k)$ is even and $k$ divides $T(k)$.
We consider the conditions $k \divides m$ and $(k+1) \divides (m+1)$.  These are
satisfied if $m \congruent k \bmod k(k+1)$, that is, $m = k + k(k+1)t$ for some $t$.
We have $m(m+1) = k(1+(k+1)t) \cdot (k+1)(1+kt)$, so that if $n = k(k+1)(1+kt)/2$ we have
$n \divides T(m)$.  Now consider $n+1 = T(k)+1 + kT(k)t$.
We have $k \divides T(k)$, so $T(k)+1$ is coprime to both $k$ and $T(k)$.
Thus the arithmetic progression $T(k)+1 + kT(k)t$ has initial term coprime to its increment
and by Dirichlet's Theorem contains infinitely many primes.
We find that there are thus infinitely many values of $t$ for which $n+1$ is prime and so
$Z(n) \le m = k + k(k+1)t$ and $Z(n+1) = n = T(k)(1+kt)$.  Hence
$$
\frac{Z(n+1)}{Z(n)} \ge \frac{n}{m} = \frac{T(k)+kT(k)t}{k+2T(k)t} > \frac{k}{3} .
$$
A similar argument holds if we consider the arithmetic progression $T(k)-1 + kT(k)t$.  
We then find infinitely many values of $t$ for which $n-1$ is prime and
$$
\frac{Z(n-1)}{Z(n)} \ge \frac{n-2}{m} = \frac{T(k)-2+kT(k)t}{k + 2T(k)t} > \frac{k}{4} .
$$
The Theorem follows by taking $k > 4L$.
\end{proof}

We note that this Theorem, combined with Lemma 1(2), gives another proof of the result
that the difference of consecutive values is unbounded.

\section{Divisibility of the pseudo-Smarandache function}

\begin{theorem}
For any integer $k \ge 2$, the equation $n/Z(n) = k$ has infinitely many solutions $n$.
\end{theorem}

\begin{proof}
Fix an integer $k \ge 2$.
Let $p$ be a prime $\congruent -1 \bmod 2k$ and put $p+1 = 2kt$.
Put $n = T(p)/t = p(p+1)/2t = p k$.  Then $n \divides T(p)$ so that $Z(n) \le p$.
We have $p \divides n$, so $Z(n) \ge p-1$: that is, $Z(n)$ must be either $p$ or $p-1$.
Suppose, if possible, that it is the latter.  In this case we have 
$2n \divides p(p+1)$ and $2n \divides (p-1)p$, so $2n$ divides $p(p+1) - (p-1)p = 2p$: 
but this is impossible since $k>1$ and so $n > p$.  We conclude that $Z(n) = p$ 
and $n/Z(n) = k$ as required.  
Further, for any given value of $k$ there are infinitely many prime values of $p$ 
satisfying the congruence condition and hence infinitely many values of $n=T(p)$ 
such $z/Z(n)=k$.
\end{proof}

\section{Another divisibility question}

\begin{theorem}
The ratio $Z(2n)/Z(n)$ is not bounded above.
\end{theorem}

\begin{proof}
Fix an integer $k$.
Let $p \congruent -1 \bmod 2^k$ be prime and put $n = T(p)$.  Then $Z(n) = p$.
Consider $Z(2n) = m$.  We have $2^k p \divides p(p+1) = 2n$ and this divides $m(m+1)/2$.
We have $m \congruent \epsilon \bmod p$ and $m \congruent \delta \bmod 2^{k+1}$
where each of $\epsilon,\delta$ can be either $0$ or $-1$.

Let $m = pt + \epsilon$.
Then $m \congruent \epsilon - t \congruent \delta \bmod 2^k$: that
is, $t \congruent \epsilon - \delta \bmod 2^k$.
This implies that either $t = 1$ or $t \ge 2^k-1$.
Now if $t=1$ then $m \le p$ and $T(m) \le T(p) = n$, which is impossible since $2n \le T(m)$.
Hence $t \ge 2^k-1$.  Since $Z(2n)/Z(n) = m/p > t/2$, we see that the ratio $Z(2n)/Z(n)$
can be made as large as desired.
\end{proof}

\section{Convergence of a series}

Ashbacher observes that the series $\sum_n 1/Z(n)^\alpha$ diverges for $\alpha=1$
and asks whether it converges for $\alpha=2$.  

In this section we prove convergence for all $\alpha > 1$.

\begin{lemma}
$$
\log n \le \sum_{m=1}^n \frac{1}{m} \le 1 + \log n ;
$$
$$
\frac12 (\log n)^2 - 0.257 \le \sum_{m=1}^n \frac{\log m}{m} \le \frac12 (\log n)^2 + 0.110
\mbox{~~for~~} n \ge 4 .
$$
\end{lemma}
\begin{proof}
For the first part, we
have $1/m \le 1/t \le 1/(m-1)$ for $t \in [m-1,m]$.  Integrating,
$$
\frac{1}{m} \le \int_{m-1}^m \frac{1}{t} \d t \le \frac{1}{m-1} .
$$
Summing,
$$
\sum_2^n \frac{1}{m} \le \int_1^n \frac{1}{t}\d t \le \sum_2^n \frac{1}{m-1} ,
$$
that is,
$$
\sum_1^n \frac{1}{m} \le 1 + \log n \mbox{~~and~~} \log n \le \sum_1^{n-1} \frac{1}{m} .
$$
The result follows.

For the second part, we similarly have $\log m/m \le \log t/t \le \log(m-1)/(m-1)$
for $t \in [m-1,m]$ when $m\ge 4$, since $\log x / x$ is monotonic decreasing for $x > \e$.
Integrating,
$$
\frac{\log m}{m} \le \int_{m-1}^m \frac{\log t}{t} \d t \le \frac{\log(m-1)}{m-1} .
$$
Summing,
$$
\sum_4^n \frac{\log m}{m} \le \int_3^n \frac{\log t}{t}\d t \le \sum_4^n \frac{\log(m-1)}{m-1} ,
$$
that is,
\begin{eqnarray*}
& &\sum_1^n \frac{\log m}{m} - \frac{\log 2}{2} - \frac{\log 3}{3} 	\\
&\le& \frac{1}{2}(\log n)^2 - \frac12(\log 3)^2 		 	\\
&\le& \sum_1^n \frac{\log m}{m} - \frac{\log n}{n} - \frac{\log 2}{2} .
\end{eqnarray*}
We approximate the numerical values
$$
\frac{\log 2}{2} + \frac{\log 3}{3} - \frac12(\log 3)^2  < 0.110
$$
and
$$
\frac{\log 2}{2} - \frac12(\log 3)^2  > -0.257 .
$$
to obtain the result.
\end{proof}

\begin{lemma}
Let $d(m)$ be the function which counts the divisors of $m$.  
For $n \ge 2$ we have
$$
\sum_{m=1}^n  d(m)/m < 7(\log n)^2 .
$$
\end{lemma}
\begin{proof}
We verify the assertion numerically for $n \le 6$.  
Now assume that $n \ge 8 > \e^2$.  We have
\begin{eqnarray*}
\sum_{m=1}^n  \frac{d(m)}{m} &=& \sum_{m=1}^n  \sum_{de=m} \frac{1}{m}
 = \sum_{d\le n} \sum_{de\le n} \frac{1}{de} 		\\
&=& \sum_{d\le n} \frac{1}{d} \sum_{e<n/d} \frac{1}{e}
   \le \sum_{d\le n} \frac{1}{d} (1+\log(n/d)) 		\\
&\le& (1 + \log n)^2 - \frac12(\log n)^2 + 0.257 	\\
&=& 1.257 + 2\log n + \frac12(\log n)^2			\\
&<& \frac43\paren(\frac{\log n}{2})^2 + 2\log n\paren(\frac{\log n}{2}) + \frac12(\log n)^2	\\
&<& 2(\log n)^2 .
\end{eqnarray*}
\end{proof}

\begin{lemma}
Fix an integer $t \ge 5$.  Let $\e^t > Y > \e^{(t-1)/2}$.
The number of integers $n$ with $\e^{t-1} < n \le \e^t$ such that $Z(n) \le Y$ is
at most $196 Y t^2$.
\end{lemma}
\begin{proof}
Consider such an $n$ with $m = Z(n) \le Y$.
Now $n \divides m(m+1)$, say $k_1n_1 = m$ and $k_2n_2 = m+1$, with $n = n_1n_2$.
Thus $k = k_1k_2 = m(m+1)/n$ and $k_1n_1 \le Y$.  The value of $k$ is bounded 
below by 2 and above by $m(m+1)/n \le 2Y^2 / \e^{t-1} = K$, say.
Given a pair $(k_1,k_2)$, the possible values of $n_1$ are bounded above by $Y/k_1$
and must satisfy the congruence condition $k_1n_1 + 1 \congruent 0$ modulo $k_2$:
there are therefore at most $Y / k_1k_2 + 1$ such values.
Since $Y/k \ge Y/K = \e^{t-1}/2Y > 1/2\e$, we have $Y/k + 1 < (2\e+1)Y/k < 7Y/k$.
Given values for $k_1,k_2$ and $n_1$, the value of $n_2$ is fixed as $n_2 = (k_1n_1+1)/k_2$.
There are thus at most $\sum_{k\le K} d(k)$ possible pairs $(k_1,k_2)$ and hence
at most $\sum_{k\le K} 7 Y d(k)/k$ possible quadruples $(k_1,k_2,n_1,n_2)$.
We have $K > 2$ so that the previous Lemma applies and we can deduce that
the number of values of $n$ satisfying the given conditions is at most $49 Y (\log K)^2$.
Now $K = 2Y^2/\e^{t-1} < 2\e^{t+1}$ so $\log K < t + 1+\log 2 < 2t$.
This establishes the claimed upper bound of $196 Y t^2$.
\end{proof}

\begin{theorem}
Fix $\frac12 < \beta < 1$ and an integer $t \ge 5$.
The number of integers $n$ with $\e^{t-1} < n \le \e^t$ such that $Z(n) < n^\beta$ is
at most $196 t^2 \e^{\beta t}$.
\end{theorem}
\begin{proof}
We apply the previous result with $Y = \e^{\beta t}$.  The conditions of $\beta$ ensure
that the previous lemma is applicable and the upper bound on the number of such $n$ is 
$ 196 \e^{\beta t} t^2 $ as claimed.
\end{proof}

\begin{theorem}
The series
$$
\sum_{n=1}^\infty \frac{1}{Z(n)^\alpha}
$$
is convergent for any $\alpha > \sqrt2$.
\end{theorem}

\begin{proof}
We note that if $\alpha > 2$ then $1/Z(n)^\alpha < 1/n^{\alpha/2}$ and the series is
convergent.  So we may assume $\sqrt2 < \alpha \le 2$.
Fix $\beta$ with $1/\alpha < \beta < \alpha/2$.
We have $\frac12 < \beta < \sqrt{\frac12} < \alpha/2$.

We split the positive integers $n>\e^4$ into two classes $A$ and $B$.  We let class $A$ be
the union of the $A_t$ where, for positive integer $t\ge5$ we put into class $A_t$
those integers $n$ such that $\e^{t-1} < n \le \e^t$ for integer $t$ and
$Z(n) \le n^\beta$.  
All values of $n$ with $Z(n) > n^\beta$ we put into class $B$.
We consider the sum of $1/Z(n)^\alpha$ over each of the two classes.
Since all terms are positive, it is sufficient to prove that each series separately
is convergent.

Firstly we observe that for $n \in B$, we have $1/Z(n)^\alpha < 1/n^{\alpha\beta}$ and
since $\alpha\beta > 1$ the series summed over the class $B$ is convergent.

Consider the elements $n$ of $A_t$: so for such $n$ we have
$\e^{t-1} < n \le \e^t$ and $Z(n) < n^\beta$.  
By the previous result, the number of values of $n$ satisfying these 
conditions is at most $196 t^2 \e^{\beta t}$.
For $n \in A_t$, we have $Z(n) \ge \sqrt{n}$,
so $1/Z(n)^\alpha \le 1/n^{\alpha/2} < 1/\e^{\alpha(t-1)/2}$.
Hence the sum of the subseries $\sum_{n \in A_t} 1/Z(n)^\alpha$
is at most $196 \e^{\alpha/2} t^2 \e^{(\beta-\alpha/2)t}$.
Since $\beta < \alpha/2$ for $\alpha > \sqrt2$, the sum over all $t$ of these terms is finite.

We conclude that $\sum_{n=1}^\infty 1/Z(n)^\alpha$ is convergent for $\alpha > \sqrt 2$
\end{proof}

\begin{theorem}
The series
$$
\sum_{n=1}^\infty \frac{1}{Z(n)^\alpha}
$$
is convergent for any $\alpha > 1$.
\end{theorem}

\begin{proof}
We fix $\beta_0 = 1 > \beta_1 > \cdots > \beta_r = \frac12$ with $\beta_j < \alpha\beta_{j+1}$ for
$0 \le j \le r-1$.
We define a partition of the integers $\e^{t-1} < n < \e^t$ into classes $B_t$ and $C_t(j)$,
$1 \le j \le r-1$.  
Into $B_t$ place those $n$ with $Z(n) > n^{\beta_1}$.  
Into $C_t(j)$ place those $n$ with $n^{\beta_{j+1}} < Z(n) < n^{\beta_j} $.  
Since $\beta_r = \frac12$ we see that every $n$ with $\e^{t-1} < n < \e^t$ is placed into one
of the classes.

The number of elements in $C_t(j)$ is at most $196 t^2\e^{\beta_j t}$ and so 
$$
\sum_{n \in C_t(j)} \frac{1}{Z(n)^\alpha} < 196 t^2 \e^{\beta_j t} \e^{-\beta_{j+1}\alpha (t-1)}
= 196 \e^{\beta_{j+1}\alpha} t^2 \e^{(\beta_j-\alpha\beta_{j+1})t} .
$$
For each $j$ we have $\beta_j < \alpha\beta_{j+1}$ so each sum over $t$ converges.

The sum over the union of the $B_t$ is bounded above by
$$
\sum_n \frac{1}{n^{\alpha\beta_1}} ,
$$
which is convergent since $\alpha\beta_1 > \beta_0 = 1$.

We conclude that $\sum_{n=1}^\infty 1/Z(n)^\alpha$ is convergent.
\end{proof}


\providecommand{\bysame}{\leavevmode\hbox to3em{\hrulefill}\thinspace}

\end{document}